# Analytical inverse for the symmetric circulant tridiagonal matrix

**Seyyed Mostafa Mousavi Janbeh Sarayi[1,]∗, Saman Tavana [2][2], Morad Karimpour [3][1], Mansour Nikkhah Bahrami4 [1*]**

[1] department of Mechanical Engineering, University of Tehran, Tehran, Iran

[2] department of Mechanical Engineering, Imperial College London, England

**Abstract** Finding the inverse of a matrix is an open problem especially when it comes to engineering problems due to their complexity and running time (cost) of matrix inversion algorithms. An optimum strategy to invert a matrix is, first, to reduce the matrix to a simple form, only then beginning a mathematical procedure. For symmetric matrices, the preferred simple form is tridiagonal. This makes tridiagonal matrices be of high interest in applied mathematics and engineering problems. This study presents a time efficient, exact analytical approach for finding the inverse, decomposition, and solving linear systems of equations where symmetric circulant matrix appears. This matrix appears in many researches and it is different from regular tridiagonal matrices as there are two corner elements. For finding the inverse matrix, a set of matrices are introduced that any symmetric circulant matrix could be decomposed into them. After that, the exact analytical inverse of this set is found which gives the inverse of circulant matrix. Moreover, solving linear equations can be carried out using implemented decomposition where this matrix appears as coefficient matrix. The method's principal strength is that it is as stable as any other direct methods (i.e., execute in a predictable number of operations). It is straightforward, understandable, solid as a rock and an exceptionally good "psychological" backup for those times when something is going wrong and you think it might be your linear-equation solver. The downside of present method, and every other direct methods, is accumulation of round off errors.

---

∗ Corresponding Author:
*Email Address*: m.mousavi.6789@gmail.com

## 1. Introduction

A $n \times n$ matrix $A = [a_{i,j}; i, j = 0,1,....,n-1]$ is said to be symmetric if $a_{i,j} = a_{j,i}$ and it is said to be tridiagonal if it has nonzero elements only on the diagonal plus or minus one column. On the other hand, a circulant matrix is a special type of Toeplitz matrix. A $n \times n$ matrix $A = [a_{i,j}; i, j = 0,1,....,n-1]$ is said to be Toeplitz if $a_{i,j} = a_{i+1,j+1}$. A Toeplitz matrix is said to be circulant if subscripts are taken modulo $n$. Thus a circulant matrix can be written as

$$A = \begin{bmatrix} a_0 & a_1 & a_2 & \cdots & a_{n-1} \\ a_{n-1} & a_0 & a_1 & \cdots & a_{n-2} \\ a_{n-2} & a_{n-1} & a_0 & \cdots & a_{n-3} \\ \vdots & \vdots & \vdots & \ddots & \vdots \\ a_1 & a_2 & a_3 & \cdots & a_0 \end{bmatrix} \quad (1)$$

In cases where only $a_0$, $a_1$ and $a_{n-1}$ are nonzero, the matrix denoted by *A* is circulant tridiagonal and if $a_1 = a_{n-1}$ the matrix is symmetric circulant.

Symmetric linear systems are prevalently appear in literature, such as electromagnetic scattering problem(Bao and Sun, 2005), molecular scattering (Poirier, 2000), structural dynamics (Feriani et al., 2000) and quantum mechanics(van Dijk and Toyama, 2007). Besides the solution of symmetric equations, (Nikkhah Bahrami and Harsini, 2007a) studied the solution of linear Toeplitz equation and then presented a solution of the linear Vandermonde equation based on the matrix decomposition method in another study (Nikkhah Bahrami and Harsini, 2007b). (Rojo, 1990) proposed an algorithm that can solve tridiagonal systems in $O(n)$ by decomposing a Toeplitz matrix into a sum of a pair of matrices, so that one of them has a simple LU-type Toeplitz like decomposition. (Vidal and Alonso, 2012) extended Rojo's algorithm to the case of symmetric Toeplitz tridiagonal equations. In the same manner, (Nemani and Garey, 2002) presented a new stable algorithm for the solution of tridiagonal circulant linear systems of equations. (Garey and Shaw, 2001) studied nonsymmetric Toepliz systems and nonsymmetric circulant systems. By making the transformation matrices, (Zheng and Shon, 2015) studied the determinant and inverse of a generalized Lucas skew circulant matrix. In another study, (Broughton and Leader, 2014) represented an analytical formula for the inverse of a symmetric circulant tridiagonal matrix as a product of a circulant matrix and its transpose.

The current study aims to present a new and comprehensive approach to decompose, calculate the inverse and solve linear systems of equations where symmetric tridiagonal and symmetric circulant tridiagonal nonsingular matrices appear. Symmetric matrices could be transformed to these tridiagonal forms using Givens rotation or householder (Press et al., 2007) method or other computational and analytical algorithms. Using present method, tridiagonal symmetric matrix could be solved for inverting problem with one step less than procedure needed for the circulant symmetric matrix. Using present method, the decomposition of symmetric circulant matrices may be found efficiently in $O(n^2)$ and the inverse of that in $O(n^{2.3728639})$ operations.

## 2. Decomposition of a symmetric circulant non-diagonal matrix

In the Eq. (1), when only $a_0$, $a_1$ and $a_{n-1}$ are nonzero, matrix $A$ is a circulant tridiagonal and if $a_1 = a_{n-1}$ the matrix is symmetric. Here we consider $a_0 = c$ and $a_1 = a_{n-1} = a$ for convenience

$$A_{n \times n} = \begin{bmatrix} c & a & 0 & \cdots & 0 & a \\ a & c & a & \ddots & \vdots & 0 \\ 0 & a & \ddots & \ddots & 0 & \vdots \\ \vdots & 0 & \ddots & \ddots & a & 0 \\ 0 & \vdots & \ddots & a & c & a \\ a & 0 & \cdots & 0 & a & c \end{bmatrix}_{n \times n}, |c| > 2|a|, a \neq 0 \qquad (2)$$

where $|c| > 2|a|$ is the diagonal dominance assumption. The normalized form of the matrix is:

$$A_{n \times n} = a \begin{bmatrix} d & 1 & 0 & \cdots & 0 & 1 \\ 1 & d & 1 & \ddots & \vdots & 0 \\ 0 & 1 & \ddots & \ddots & 0 & \vdots \\ \vdots & 0 & 1 & \ddots & 1 & 0 \\ 0 & \vdots & \ddots & 1 & d & 1 \\ 1 & 0 & \cdots & 0 & 1 & d \end{bmatrix}_{n \times n} = a \times \overline{A_{n \times n}} \qquad (3)$$

where $d = c/a$. Let matrix $K_{n \times n}$ be of the form

$$K_{n \times n} = \begin{bmatrix} f_1(d) & 0 & \cdots & \cdots & 0 \\ f_1(d) & f_2(d) & \ddots & \ddots & \vdots \\ \vdots & f_2(d) & \ddots & \ddots & \vdots \\ \vdots & \vdots & \ddots & f_{n-1}(d) & 0 \\ f_1(d) & f_2(d) & \cdots & f_{n-1}(d) & f_n(d) \end{bmatrix}_{n \times n} \qquad (4)$$

Where matrix $K$ is a lower triangular alternant matrix and $f_i$ are obtained from the following recurrence relation:

$$f_{i+1}(d) = -df_i(d) - f_{i-1}(d), \\ f_0(d) = 0, f_1(d) = 1, i = 1,2,3,\ldots,n-1, f_{i+1}(d) \neq 0 \qquad (5)$$

The bidiagonal inverse of matrix $K$ could be obtained as:

$$K^{-1}_{n\times n} = \begin{bmatrix} \dfrac{1}{f_1(d)} & 0 & \cdots & \cdots & 0 \\ -\dfrac{1}{f_2(d)} & \dfrac{1}{f_2(d)} & \ddots & \ddots & \vdots \\ 0 & -\dfrac{1}{f_3(d)} & \ddots & \ddots & \vdots \\ \vdots & \ddots & \ddots & \dfrac{1}{f_{n-1}(d)} & 0 \\ 0 & \cdots & 0 & -\dfrac{1}{f_n(d)} & \dfrac{1}{f_n(d)} \end{bmatrix}_{n\times n} \quad (6)$$

Based on Eq. (6), all entries of matrix $K^{-1}_{n\times n}$ are zero except for the elements of the main diagonal and subdiagonal. It can be seen in Eqs. (4) and (6) that matrix $K$ and its inverse have elements that can be calculated easily using Eq. (5).

Let matrix $R_{n\times n}$ be of the form

$$R_{n\times n} = \begin{bmatrix} 1 & 0 & \cdots & 0 & 0 \\ 0 & 1 & \ddots & \vdots & 0 \\ \vdots & \ddots & \ddots & 0 & \vdots \\ 0 & \cdots & 0 & 1 & 0 \\ r_1(d) & r_2(d) & \cdots & r_{n-1}(d) & 1 \end{bmatrix}_{n\times n} \quad (7)$$

where matrix $R$ is a sparse matrix in a way that the entries of main diagonal are one and the entries of its last row are $r_j$ except for the entry of the $n^{th}$ column which is one. In Eq. (7) $r_j$ are:

$$r_j(d) = \dfrac{f_n(d)f_1(d)}{f_{j+1}(d)f_j(d)}, \quad j=1,2,3,\ldots,n-1 \quad (8)$$

The inverse of the matrix $R$ can be calculated using decomposition method:

$$R^{-1}_{n\times n} = \begin{bmatrix} 1 & 0 & \cdots & 0 & 0 \\ 0 & 1 & \ddots & \vdots & 0 \\ \vdots & \ddots & \ddots & 0 & \vdots \\ 0 & \cdots & 0 & 1 & 0 \\ -r_1(d) & -r_2(d) & \cdots & -r_{n-1}(d) & 1 \end{bmatrix}. \quad (9)$$

Now let introduce matrix $A_1^{kR}$

$$A_1^{kR} = \begin{bmatrix} -f_2(d) & 0 & 0 & \cdots & 0 & 1 \\ f_1(d) & -f_3(d) & 0 & \cdots & 0 & 1 \\ 0 & f_2(d) & -f_4(d) & \ddots & \vdots & \vdots \\ \vdots & & \ddots & \ddots & 0 & 1 \\ \vdots & \vdots & & \ddots & -f_n(d) & 1 \\ 0 & \cdots & & \cdots & 0 & f_{n-1}(d)+1 & g_{n+1}(d) \end{bmatrix}_{n \times n} \quad (10)$$

In this matrix only one element should be computed which is $g_{n+1}(d)$ to generate matrix $A_1^{kR}{}_{n \times n}$ and it can be obtained as follows:

$$g_{n+1}(d) = 1 - f_{n+1}(d) + \sum_{j=1}^{n-1} r_j(d) + r_{n-1}(d) f_{n-1}(d). \quad (11)$$

or

$$g_{n+1}(d) = 1 + f_1(d) - f_{n+1}(d) + \sum_{j=1}^{n-1} r_j(d), f_1(d) = 1. \quad (12)$$

Matrix $A_1^{kR}$ is a lower triangular matrix and its inverse is

$$(A_1^{kR})^{-1} = \begin{bmatrix} \dfrac{1}{f_2(d)} & 0 & 0 & 0 & \cdots & 0 & 0 \\ \dfrac{f_1(d)}{f_2(d)f_3(d)} & \dfrac{1}{f_3(d)} & 0 & 0 & \cdots & 0 & 0 \\ \dfrac{f_1(d)}{f_3(d)f_4(d)} & \dfrac{f_2(d)}{f_3(d)f_4(d)} & \dfrac{1}{f_4(d)} & 0 & \ddots & \vdots & \vdots \\ & & & \ddots & \ddots & 0 & 0 \\ \vdots & \vdots & & \ddots & \ddots & \ddots & 0 \\ \dfrac{f_1(d)}{f_{n-1}(d)f_n(d)} & \dfrac{f_2(d)}{f_{n-1}(d)f_n(d)} & \cdots & \dfrac{f_{n-2}(d)}{f_{n-1}(d)f_n(d)} & \dfrac{1}{f_n(d)} & 0 \\ \dfrac{f_1(d)}{f_n(d)g_{n+1}(d)} & \dfrac{f_2(d)}{f_n(d)g_{n+1}(d)} & \cdots & \dfrac{f_{n-2}(d)}{f_n(d)g_{n+1}(d)} & \dfrac{f_{n-1}(d)}{f_n(d)g_{n+1}(d)} & \dfrac{1}{g_{n+1}(d)} \end{bmatrix}_{n \times n} \quad (13)$$

In order to illustrate the decomposition of matrix $A$ to these three matrices, we have:

$$A = a\ K^{-1} R^{-1} \left(A_1^{kR}\right)^T \quad (14)$$

Eq. (14) is a decomposition of matrix A into three sparse and easy to calculate matrices. $K^{-1}$ is presented in Eq. (6), $R^{-1}$ in Eq. (9), and $A_1^{KR}$ using Eq. (10). Using this decomposition, the inverse of matrix A can easily be derived as follow:

$$A^{-1} = \frac{1}{a}\left[\left(A_1^{kR}\right)^{-1}\right]^T RK \qquad (15)$$

However, based on Eq. (15), linear system of $AX = b$ may be solved as $X = A^{-1}b$. In addition, as a different approach to solve $AX = b$, we have the following relation based on Eq. (14):

$$AX = b \rightarrow aRK\bar{A}X = \bar{b} \rightarrow \left(A_1^{kR}\right)^T X = \bar{b} \rightarrow A^{kR}X = \bar{b} \qquad (16)$$

where $\bar{b} = aRKb$. With respect to the last equation of Eq. (16) which is an updated form of $AX = b$, the values of $X$ can be easily obtained using gauss method (back substitution).

## 3. Numerical Example

Now let us consider the following 5x5 ciculant tridiagonal symmetric matrix:

$$A = \begin{bmatrix} 5 & 2 & 0 & 0 & 2 \\ 2 & 5 & 2 & 0 & 0 \\ 0 & 2 & 5 & 2 & 0 \\ 0 & 0 & 2 & 5 & 2 \\ 2 & 0 & 0 & 2 & 5 \end{bmatrix}, A = \begin{bmatrix} 5 & 2 & 0 & 0 & 2 \\ 2 & 5 & 2 & 0 & 0 \\ 0 & 2 & 5 & 2 & 0 \\ 0 & 0 & 2 & 5 & 2 \\ 2 & 0 & 0 & 2 & 5 \end{bmatrix}, \bar{A} = \begin{bmatrix} 2.5 & 1 & 0 & 0 & 1 \\ 1 & 2.5 & 1 & 0 & 0 \\ 0 & 1 & 2.5 & 1 & 0 \\ 0 & 0 & 1 & 2.5 & 1 \\ 1 & 0 & 0 & 1 & 2.5 \end{bmatrix}, d = 2.5, a = 2$$

From Eq. (5) $f_i$ are calculated as follow:

$$f_2 = -2.5,\ f_3 = 5.25,\ f_4 = -10.625,\ f_5 = 21.3125,\ f_6 = -42.6566$$

Subsequently matrix $K$ is calculated by using Eq. (4):

$$K = \begin{bmatrix} 1 & 0 & 0 & 0 & 0 \\ 1 & -2.5 & 0 & 0 & 0 \\ 1 & -2.5 & 5.25 & 0 & 0 \\ 1 & -2.5 & 5.25 & -10.625 & 0 \\ 1 & -2.5 & 5.25 & -10.625 & 21.3125 \end{bmatrix}$$

Inverse of this matrix can easily be obtained by using Eq. (6):

$$K^{-1} = \begin{bmatrix} 1 & 0 & 0 & 0 & 0 \\ 0.4 & -0.4 & 0 & 0 & 0 \\ 0 & -0.1905 & 0.1905 & 0 & 0 \\ 0 & 0 & 0.0941 & -0.0941 & 0 \\ 0 & 0 & 0 & -0.0469 & 0.0469 \end{bmatrix}$$

Matrix $R$ from Eq. (7):

$$R = \begin{bmatrix} 1 & 0 & 0 & 0 & 0 \\ 0 & 1 & 0 & 0 & 0 \\ 0 & 0 & 1 & 0 & 0 \\ 0 & 0 & 0 & 1 & 0 \\ -8.525 & -1.6238 & -0.3821 & -0.0941 & 1 \end{bmatrix}$$

Considering Eq. (9) its inverse is

$$R^{-1} = \begin{bmatrix} 1 & 0 & 0 & 0 & 0 \\ 0 & 1 & 0 & 0 & 0 \\ 0 & 0 & 1 & 0 & 0 \\ 0 & 0 & 0 & 1 & 0 \\ 8.525 & 1.6238 & 0.3821 & 0.0941 & 1 \end{bmatrix}$$

Using Eq. (10), matrix $A_1^{kR}$ is obtained ($g_6(d) = 68.0625$):

$$A_1^{KR} = \begin{bmatrix} 5 & 0 & 0 & 0 & 0 \\ 2 & -10.5 & 0 & 0 & 0 \\ 0 & -5 & 21.25 & 0 & 0 \\ 0 & 0 & 10.5 & -42.625 & 0 \\ 2 & 2 & 2 & -19.25 & 68.0625 \end{bmatrix}$$

Finally inverse of matrix $A$ is calculated using Eq. (15):

$$A^{-1} = \begin{bmatrix} 0.3131 & -0.1414 & 0.0404 & 0.0404 & -0.1414 \\ -0.1414 & 0.3131 & -0.1414 & 0.0404 & 0.0404 \\ 0.1414 & -0.1414 & 0.3131 & -0.1414 & 0.0404 \\ 0.1414 & 0.0404 & -0.1414 & 0.3131 & -0.1414 \\ -0.1414 & 0.0404 & 0.0404 & -0.1414 & 0.3131 \end{bmatrix}$$

## 4. Conclusion

This study presents a new comprehensive approach for finding the inverse of symmetric circulant tridiagonal matrices. As such, three $K^{-1}$, $R^{-1}$ and $\left(A_1^{kR}\right)^T$ were introduced as component matrices and their inverses are found analytically. The suggested decomposition method for special matrices could be applied for solving related linear equations. The proposed method in this study is even simpler in some cases including tridiagonal matrices (not circulant tridiagonal) since there is no need for computing matrix R and its inverse as they are identity matrix. In this situation, only matrices $K^{-1}$, $\left(A_1^{kR}\right)^T$ and their inverses give the inverse matrix. Further studies on the other special cases like Toeplitz and sparse matrices could be done and have their decomposition and inverse matrices found.

## 5. References

Bao, G., Sun, W., 2005. A Fast Algorithm for the Electromagnetic Scattering from a Large Cavity. SIAM J. Sci. Comput. 27, 553–574. doi:10.1137/S1064827503428539

Broughton, S.A., Leader, J.J., 2014. Analytical Solution of the Symmetric Circulant Tridiagonal Linear System. Math. Sci. Tech. Reports 103. doi:http://scholar.rose-hulman.edu/math_mstr/103

Feriani, A., Perotti, F., Simoncini, V., 2000. Iterative system solvers for the frequency analysis of linear mechanical systems. Comput. Methods Appl. Mech. Eng. 190, 1719–1739. doi:10.1016/S0045-7825(00)00187-0

Garey, L.E., Shaw, R.E., 2001. A parallel method for linear equations with tridiagonal Toeplitz coefficient matrices. Comput. Math. with Appl. 42, 1–11. doi:10.1016/S0898-1221(01)00125-0

Nemani, S.S., Garey, L.E., 2002. Parallel algorithms for solving tridiagonal and near-circulant systems. Appl. Math. Comput. 130, 285–294. doi:10.1016/S0096-3003(01)00096-0

Nikkhah Bahrami, M., Harsini, I., 2007a. Solution of linear Toeplitz equation by matrix decomposition method. Pasific Conf. Comput. Math. 1, 693.

Nikkhah Bahrami, M., Harsini, I., 2007b. Solution of linear Vandermond equation by matrix

decomposition method. Asian Pasific Conf. Comput. Math. 1, 705.

Poirier, B., 2000. Efficient preconditioning scheme for block partitioned matrices with structured sparsity. Numer. Linear Algebr. with Appl. 7, 715–726. doi:10.1002/1099-1506(200010/12)7:7/8<715::AID-NLA220>3.0.CO;2-R

Press, W.H., Teukolsky, S. a, Vetterling, W.T., Flannery, B.P., 2007. Numerical Recepies: The Art of Scientific Computing, Third. ed, Cambridge University Press. Press Syndicate of the University of Cambridge. doi:http://numerical.recipes/

Rojo, O., 1990. A new method for solving symmetric circulant tridiagonal systems of linear equations. Comput. Math. with Appl. 20, 61–67. doi:10.1016/0898-1221(90)90165-G

van Dijk, W., Toyama, F.M., 2007. Accurate numerical solutions of the time-dependent Schrödinger equation. Phys. Rev. E 75, 36707. doi:10.1103/PhysRevE.75.036707

Vidal, A.M., Alonso, P., 2012. Solving systems of symmetric Toeplitz tridiagonal equations: Rojo's algorithm revisited. Appl. Math. Comput. 219, 1874–1889. doi:10.1016/j.amc.2012.08.030

Zheng, Y., Shon, S., 2015. Exact determinants and inverses of generalized Lucas skew circulant type matrices. Appl. Math. Comput. 270, 105–113. doi:10.1016/j.amc.2015.08.021